\documentclass[12pt]{amsart}
\usepackage{a4wide}
\usepackage{graphicx}
\usepackage{epstopdf}
\usepackage{amssymb,amsmath, enumerate}
\usepackage[pdfborder={0 0 0}]{hyperref}

\newtheorem{pro}{Proposition}[section]
\theoremstyle{definition}

\newtheorem{rem}{Remark}[section]
\title{Can an infinite left-product of nonnegative matrices be expressed in terms of infinite left-products of stochastic ones?}
\author[A. Thomas]{Alain Thomas}
\address[Alain Thomas]{
LATP, 39, rue Joliot-Curie,\hfil\break
13453 Marseille, Cedex 13,
France}
\email{thomas@cmi.univ-mrs.fr}
\date{}  
\begin{document}
\baselineskip=18pt
\maketitle
\begin{abstract}
If a left-product $M_n\dots M_1$ of square complex matrices converges to a nonnull limit when $n\to\infty$ and if the $M_n$ belong to a finite set, it is clear that there exists an integer $n_0$ such that the $M_n$, $n\ge n_0$, have a common right-eigenvector $V$ for the eigenvalue~$1$. Now suppose that the $M_n$ are nonnegative and that $V$ has positive entries. Denoting by $\Delta$ the diagonal matrix whose diagonal entries are the entries of $V$, the stochastic matrices $S_n=\Delta^{-1}M_n\Delta$ satisfy $M_n\dots M_{n_0}=\Delta S_n\dots S_{n_0}\Delta^{-1}$, so the problem of the convergence of $M_n\dots M_1$ reduces to the one of $S_n\dots S_{n_0}$. In this paper we still suppose that the $M_n$ are nonnegative but we do not suppose that $V$ has positive entries. The first section details the case of the $2\times2$ matrices, and the last gives a first approach in the case of $d\times d$ matrices.
\end{abstract}
\par\par

\vskip20pt

\centerline{\sc Introduction}

\smallskip

The problem of the convergence of $M_1\dots M_n$ or $M_n\dots M_1$ for all the sequences $(M_n)$ with terms in a finite set of complex matrices, is studied for instance in \cite{DL}, \cite{DL'} and \cite{EF}. The same problem in case of stochastic matrices is also classical, see for instance \cite[chapter 4]{S}.

On the other side there exist much results about the distribution of the product matrix $M_1\dots M_n$ where the $M_i$ are taken in a set of stochastic matrices endowed with some probability measure. In \cite{MNR} and \cite{MN}, Mukherjea, Nakassis and Ratti give conditions for the limit distribution of $M_1\dots M_n$ to be discrete or continuous singular; this contains for instance the case of the Erd\H os measure \cite{Erd}. It is well known that in much cases the normalized product $M_1\dots M_n$, if the $M_i$ are taken in a finite set of nonnegative matrices endowed with some positive probability $P$, converges $P$-almost everywhere to a rank one matrix; but this general result should be more consistent if one can to specify the \hbox{$P$-negligeable} set of divergence.

Notice that the Erd\H os measure is studied more in detail in \cite{OST} by an other method, using a finite set of matrices and the asymptotic properties of the columns in the products of matrices taken in this set.

In the present paper we first consider (\S1) the left-products and the right-products of $2\times2$ matrices (resp. $2\times2$ stochastic matrices). For the left-products of stochastic matrices, the hypothesis that the matrices belong to a finite set is not necessary. We recover the known results by giving several formulations of the necessary and sufficient conditions of convergence.

In \S2 we associate to any sequence $(M_n)$ of nonnegative $d\times d$ matrices, some sequences of stochastic ones, let $\left(S^{(i)}_n\right)$  for $1\le i\le t$. The convergence of $M_n\dots M_1$ is equivalent to the one of the $S^{(i)}_n\dots S^{(i)}_1$ and some additional condition. 

\section{Products of $2\times2$ matrices}

\subsection{Convergent left-products}

Suppose that the left-product $M_n\dots M_1$ of some nonnegative matrices $M_n=\left(\begin{array}{cc}a_n&b_n\\ c_n&d_n\end{array}\right)$ converges to a nonnull limit, let $Q$, and that the set $\{M\;;\;\exists n, M=M_n\}$ is finite. It is clear that from a rank~$n_0$, the $M_n$ belong to the set
$$
\{M\;;\;\hbox{ there exists infinitely many }n\hbox{ such that }M=M_n\}
$$
and consequently, for any $M$ in this set, $MQ=M$. In other words the nonnull columns of $Q$ are eigenvectors -- for the eigenvalue~$1$~-- of any $M_n$, $n\ge n_0$. If for instance this eigenvector is $\left(\begin{array}{cc}v_1\\0\end{array}\right)$, then $M_n=\left(\begin{array}{cc}1&b_n\\ 0&d_n\end{array}\right)$ and $M_n\dots M_{n_0}=\left(\begin{array}{cc}1&\sum_{i=n_0}^nb_id_{i-1}\dots d_{n_0}\\ 0&d_n\dots d_{n_0}\end{array}\right)$

\begin{pro}
$M_n\dots M_{n_0}$ converges to a nonnull limit, for any fixed $n_0$ and when $n\to\infty$, if and only if

$\_$ either the $M_n$ have from a certain rank, a common right-eigenvector $\left(\begin{array}{cc}v_1\\v_2\end{array}\right)$ for the eigenvalue $1$, with $v_1v_2>0$, and the left-product of the stochastic matrices $S_n=\left(\begin{array}{cc}1\over v_1&0\\ 0&1\over v_2\end{array}\right)M_n\left(\begin{array}{cc} v_1&0\\ 0&v_2\end{array}\right)$ from any rank $n_0$ converges;

$\_$ or the $M_n$ have from a certain rank $\left(\begin{array}{cc}1\\0\end{array}\right)$ for common positive right-eigenvector with respect to the eigenvalue $1$, the sum $\displaystyle\sum_{n=n_0}^\infty b_nd_{n-1}\dots d_{n_0}$ is finite and $d_n\dots d_{n_0}$ converges, for any $n_0$, when $n\to\infty$;

$\_$ or the $M_n$ have from a certain rank $\left(\begin{array}{cc}0\\1\end{array}\right)$ for common positive right-eigenvector with respect to the eigenvalue $1$, the sum $\displaystyle\sum_{n=n_0}^\infty c_na_{n-1}\dots a_{n_0}$ is finite and $a_n\dots a_{n_0}$  converges, for any $n_0$, when $n\to\infty$.

\end{pro}

\subsection{Case of stochastic matrices}

The case of the matrices
$$
S_n=\left(\begin{array}{cc}x_n&1-x_n\\y_n&1-y_n\end{array}\right),\ x_n,y_n\in[0,1],
$$
is also trivial because one can compute the left-product
\begin{equation}\label{Q_n}
Q_n:=S_n\dots S_1=\left(\begin{array}{cc}t_n&1-t_n\\s_n&1-s_n\end{array}\right)\quad\hbox{where }\quad\left\{\begin{array}{l}s_n:=\sum_{i=1}^ny_i\det Q_{i-1}\\ t_n:=s_n+\det Q_n\\Q_0:=I.\end{array}\right.
\end{equation}
To find the conditions for the sequence of matrices $(Q_n)$ to converge one can use the relation
\begin{equation}\label{ssss}
s_n=s_{n_0-1}+s_{n_0,n}\det Q_{n_0-1},\quad\hbox{where }s_{n_0,n}:=\sum_{i=n_0}^n\left(y_i\prod_{n_0\le j<i}\det S_j\right).
\end{equation}
$s_{n_0,n}$ belongs to $[0,1]$ because it is one of the entries of the stochastic matrix $S_n\dots S_{n_0}$. Hence, in case $\det Q_n$ has limit $0$ the relation (\ref{ssss}) implies that $(s_n)$ is Cauchy; so $(s_n)$, $(t_n)$ and $(Q_n)$ converge.

Suppose now that $\det Q_n$ do not have limit $0$. Since $\det Q_n=\prod_{i=1}^n(x_i-y_i)$ with $x_i-y_i\in[-1,1]$, the non-increasing sequence $(\vert\det Q_n\vert)$ has a positive limit $\delta$ hence $\vert x_n-y_n\vert$ has limit~$1$; $(x_n,y_n)$ cannot have other limit points than $(0,1)$ and $(1,0)$.

In case $(0,1)$ is one of its limit points, $y_n\vert\det Q_{n-1}\vert$ do not tend to $0$ hence the series $\sum_ny_n\det Q_{n-1}$ diverges and $(Q_n)$ also do.

In case $(1,0)$ is the unique limit point of $(x_n,y_n)$, $x_n-y_n$ is positive from a rank $n_0$. Since the series $\sum_n\log\vert x_n-y_n\vert$ converges to $\log\delta$, the inequalities $\sum_{n\ge n_0}\log(x_n-y_n)\le\sum_{n\ge n_0}\log(1-y_n)\le-\sum_{n\ge n_0}y_n$ prove that the series $\sum_ny_n$ converges. Since $\det Q_n$ has limit $\delta$ or $-\delta$ according to the sign of $\det Q_{n_0-1}$, the sequences $(s_n)$, $(t_n)$ and $(Q_n)$ converge.
\vskip10pt
Consider now the right-product $P_n:=S_1\dots S_n$ and suppose that the $S_n$ belong to a finite set. As noticed in \cite{EF} it is clear that the nonnull rows of the limit matrix $P$, if this matrix exists, are nonnegative left-eigenvectors -- for the eigenvalue $1$ -- of each matrix $S$ such that $S_n=S$ for infinitely many $n$, because the equality $P_n=P_{n-1}S$ implies $P=PS$. So $P_n$ can converge only if the $S_n$, for $n$ greater or equal to some integer $n_0$, have a common nonnegative left-eigenvector for the eigenvalue $1$.

We suppose there exists such a left-eigenvector, let $L$, and we search the condition for $(P_n)$ to converge. Notice that the $S_n$ for $n\ge n_0$ commute: since any $2\times2$ stochastic matrix $S$ has left-eigenvector $\left(\begin{array}{cc}1&-1\end{array}\right)$ for the eigenvalue $\det S$, both vectors $L$ and $\left(\begin{array}{cc}1&-1\end{array}\right)$ are orthogonal to the columns of $S_nS_{n'}-S_{n'}S_n$ and consequently this matrix is null. So we fall again in the case of the left-products.

In case $S_{n_0}\dots S_n$ diverges, nevertheless the sequence of the row-vectors $LS_{n_0}\dots S_n$ converges (it is constant). Let $L'$ be some row-eigenvector not colinear to $L$; considering the invertible matrix $M$ whose rows are $L$ and $L'$, $MS_{n_0}\dots S_n$ obviously diverges hence $L'S_{n_0}\dots S_n$ also do. Consequently $S_1\dots S_n$ diverges if and only if at least one of the rows of $S_1\dots S_{n_0-1}$ is not colinear to $L$. 

We have proved the following

\begin{pro}Let $(S_n)$ be a sequence of $2\times2$ stochastic matrices, namely
$$
S_n=\left(\begin{array}{cc}x_n&1-x_n\\y_n&1-y_n\end{array}\right).
$$
\vskip5pt
(i) The left-product $Q_n=S_n\dots S_1$ converges if and only if $\prod_{k=1}^n(x_k-y_k)$ has limit $0$ or $(x_n,y_n)$ has limit $(1,0)$ when $n\to+\infty$.
\vskip5pt
(ii) It diverges only in the case where $\sum_n(1-\vert x_n-y_n\vert)$ converges and $S_n$ do not tend to $\left(\begin{array}{cc}1&0\\0&1\end{array}\right)$.
\vskip5pt
(iii) If it converges, $\lim_{n\to\infty}Q_n=\left(\begin{array}{cc}s+q&1-s-q\\s&1-s\end{array}\right)$ where $s:=\sum_ny_n\det Q_{n-1}$ and $q:=\lim_{n\to\infty}(\det Q_n)$.
\vskip5pt
(iv) Suppose now that the $S_n$ belong to a finite set. Then the right-product $P_n=S_1\dots S_n$ converges if and only if there exists $n_0$ such that

$\_$ the matrices $S_n$ for $n\ge n_0$ have a common left-eigenvector for the eigenvalue~$1$

$\_$ and either $\lim_{n\to\infty}\prod_{k=n_0}^n(x_k-y_k)=0$, or $\lim_{n\to\infty}(x_n,y_n)=(1,0)$, or the rows of $S_1\dots S_{n_0-1}$ are colinear to the left-eigenvector.
\vskip5pt
(v) Suppose the $S_n$, $n\ge1$, have a common left-eigenvector with respect to the eigenvalue~$1$. Then the $S_n$ commute and $P_n=Q_n$.
\end{pro}

\section{Left products of $d\times d$ nonnegative matrices}

Let us first give one example in order to illustrate the proposition that follows: we consider the products $Q_n=M_n\dots M_1$, where the $M_i$ belong to the set of nonnegative matrices of the form
$$
M=\left(\begin{array}{ccccccc}a&b&3-3a-2b&e&4-4a-3b-5e&1&0\\c&d&2-3c-2d&f&3-4c-3d-5f&0&1\\0&0&1&0&0&1&0\\0&0&0&g&5-5g&0&1\\0&0&0&h&1-5h&1&0\\0&0&0&0&0&x&1-x\\0&0&0&0&0&y&\displaystyle{1\over2}-y\end{array}\right)
$$
where $a,b,c,d,e,f,g,h,x,y$ are some reals such that $M$ has exactly twenty four nonnull entries. Since the eigenspace associated to the eigenvalue $1$ is generated by $\left(\begin{array}{c}3\\2\\1\\0\\0\\0\\0\end{array}\right)$ and $\left(\begin{array}{c}4\\3\\0\\5\\1\\0\\0\end{array}\right)$, we can associate to each matrix $M$ two submatrices with positive eigenvectors:

the submatrix $M^{\{1,2,3\}}$ of the entries of $M$ with row and column indexes in $\{1,2,3\}$ and the submatrix $M^{\{1,2,4,5\}}$ of the entries of $M$ with row and column indexes in $\{1,2,4,5\}$. Then we associate two stochastic matrices $S=\Delta^{-1}M^{\{1,2,3\}}\Delta$ and $S'={\Delta'}^{-1}M^{\{1,2,3\}}\Delta'$, where $\Delta=\left(\begin{array}{ccc}3&0&0\\0&2&0\\0&0&1\end{array}\right)$ and $\Delta'=\left(\begin{array}{cccc}4&0&0&0\\0&3&0&0\\0&0&5&0\\0&0&0&1\end{array}\right)$; namely$$
S=\left(\begin{array}{ccc}a&2b/3&1-a-2b/3\\3c/2&d&1-3c/2-d\\0&0&1\end{array}\right)\quad\hbox{and}\quad S'=\left(\begin{array}{cccc}a&3b/4&5e/4&1-a-3b/4-5e/4\\4c/3&d&5f/3&1-4c/3-d-5f/3\\0&0&g&1-g\\0&0&5h&1-5h\end{array}\right).
$$
Now we obtain $21$ of the $49$ entries of $Q_n=M_n\dots M_1$ in function of two products of stochastic matrices, and the other entries in the first five columns of $Q_n$ are null: indeed $Q_n$ has for submatrices $M^{\{1,2,3\}}_n\dots M^{\{1,2,3\}}_1=\Delta S_n\dots S_1\Delta^{-1}$ and $M^{\{1,2,4,5\}}_n\dots M^{\{1,2,4,5\}}_1=\Delta'S'_n\dots S'_1{\Delta'}^{-1}$.

The products $S_n\dots S_1$ and $S'_n\dots S'_1$ converge to rank $1$ matrices: use for instance \cite{L}, or use the previous section and the formula for the products of triangular-by-blocks matrices that is,
$$
\prod_{i=n}^1\left(\begin{array}{cc}A_i&B_i\\0&D_i\end{array}\right)=\left(\begin{array}{cc}\prod_{i=n}^1A_i&\sum_{i=1}^nA_n\dots A_{i+1}B_iD_{i-1}\dots D_1\\0&\prod_{i=n}^1D_i\end{array}\right).
$$
The limit of $Q_n$ is a rank $2$ matrix of the form $\left(\begin{array}{ccccccc}0&0&3&4\alpha&4\beta&4\gamma+6&4\delta+6\\0&0&2&3\alpha&3\beta&3\gamma+4&3\delta+4\\0&0&1&0&0&2&2\\0&0&0&5\alpha&5\beta&5\gamma&5\delta\\0&0&0&\alpha&\beta&\gamma&\delta\\0&0&0&0&0&0&0\\0&0&0&0&0&0&0\end{array}\right)$.

The following proposition generalizes what we have seen on the example. For any $d\times d$ matrix $M$ and for any $K,K'\subset\{1,\dots,d\}$ we denote by $M^K$ the submatrix of the entries of $M$ whose row index and column index belong to $K$ and by $M^{K,K'}$ the submatrix of the entries of $M$ whose row index belongs to $K$ and column index belongs to $K'$. By commodity we use the same notation for the row-matrices $L$ (resp. the column-matrices~$V$): $L^K$ (resp. $V^K$) is the row-matrice (resp. column-matrice) of the entries of $L$ (resp. $V$) whose index belong to $K$.

\begin{pro}Let $(M_n)$ be a sequence of nonnegative $d\times d$ matrices that belong to a given finite set. The left-product $Q_{n,n_0}:=M_n\dots M_{n_0}$ converges to a nonnull limit when $n\to\infty$ -- for each positive integer $n_0$ -- if and only if $(M_n)$ satisfies both conditions:
\vskip5pt
(i) there exist some subsets of $\{1,\dots,d\}$, let $K_1,\dots,K_t$ with complementaries $K_1^c,\dots,K_t^c$, and some diagonal matrices with positive diagonals, let $\Delta_1,\dots,\Delta_t$, such that the $S_n^{(i)}={\Delta_i}^{-1}M_n^{K_i}\Delta_i$ are stochastic, the $M_n^{K_i^c,K_i}$ are null, and $\lim_{n\to\infty}S_n^{(i)}\dots S_{n_0}^{(i)}$ exists for any $i$ and $n_0$;
\vskip5pt
(ii) setting $K=\cup_iK_i$, $\lim_{n\to\infty}M^{K^c}_n\dots M^{K^c}_{n_0}$ is the null matrix for any $n_0$ and the series $\sum_{i=n_0}^\infty M^K_n\dots M^K_{i+1}M^{K,K^c}_iM^{K^c}_{i-1}\dots M^{K^c}_{n_0}$ converges.

\end{pro}

\begin{proof}If the conditions (i) and (ii) are satisfied, the entries of $M_n\dots M_{n_0}$ with column index in $K$ converges either to $0$ or to the entries of the matrices $\lim_{n\to\infty}\Delta_iS_n^{(i)}\dots S_{n_0}^{(i)}{\Delta_i}^{-1}$, $i=1,\dots,t$. Consequently $M_n\dots M_{n_0}$ converges, by using the formula of product of triangular-by-blocs matrices:
$$
\prod_{i=n}^{n_0}\left(\begin{array}{cc}A_i&B_i\\0&D_i\end{array}\right)=\left(\begin{array}{cc}\prod_{i=n}^{n_0}A_i&\sum_{i=n_0}^nA_n\dots A_{i+1}B_iD_{i-1}\dots D_{n_0}\\0&\prod_{i=n}^{n_0}D_i\end{array}\right).
$$

Conversely suppose that $Q_{n,n_0}=M_n\dots M_{n_0}$ converges to a nonnull limit for any $n_0\in\mathbb N$. Since the $M_n$ belong to a finite set, we can choose $n_0$ large enough such that all the matrices $M$ that are equal to $M_n$ for at least one $n\ge n_0$, are also equal to $M_n$ for infinitely many $n$. Then the nonnull columns of the limit matrix $Q$ are right-eigenvectors of all the $M_n$, $n\ge n_0$, for the eigenvalue $1$, because the equality $Q_{n,n_0}=M_nQ_{n-1,n_0}=MQ_{n-1,n_0}$ implies $Q=MQ$. Let $t$ be the rank of $Q$, $t=0$ if $Q$ is null, and denote by $V_1,\dots,V_t$ the linearly independent columns of $Q$. We denote also by $K_i$ the set of the indexes of the nonnull entries in $V_i$, and by $\Delta_i$ the diagonal matrix whose diagonal entries are the nonnull entries of $V_i$. Now (i) results from the equality $M_n^{D,K_i}V_i^{K_i}=V_i$, $D:=\{1,\dots,d\}$; $S_n^{(i)}\dots S_{n_0}^{(i)}$ converges because $\Delta_iS_n^{(i)}\dots S_{n_0}^{(i)}{\Delta_i}^{-1}$ is a submatrix of $Q_{n,n_0}$.

On the other side, denoting by $K$ the union of the $K_i$, the product $M^{K^c}_n\dots M^{K^c}_{n_0}$ and the sum of products $\sum_{i=n_0}^nM^K_n\dots M^K_{i+1}M^{K,K^c}_iM^{K^c}_{i-1}\dots M^{K^c}_{n_0}$ converge when $n\to\infty$ because they are submatrices of $Q_{n,n_0}$. The first converges to $0$: by the definitions of $K$ and the vectors $V_i$, the rows of Q whose indexes belong to $K^c$ are null.

\end{proof}

\begin{rem}The condition that $\sum_{i=1}^\infty M^K_n\dots M^K_{i+1}M^{K,K^c}_iM^{K^c}_{i-1}\dots M^{K^c}_1$ converges cannot be avoided. Suppose for instance that $M^K_n\dots M^K_1$ converges and that $M^{K,K^c}_n$ is for any $n$ the identity ${d\over2}\times{d\over2}$ matrix, let $I_{d\over2}$, $d$ even. Suppose also $M^{K^c}_n=d_nI_{d\over2}$ for any $n$, where the positive reals $d_n$ satisfy $\lim_{n\to\infty}d_n\dots d_1=0$ and \hbox{$\sum_{n=1}^\infty d_n\dots d_1=\infty$}. Then if the $M_n$ have a common right-eigenvector for the eigenvalue $1$, it has the form $\begin{pmatrix}W\\\theta_{d\over2}\end{pmatrix}$ where $W$ is an eigenvector of the $M_n^K$ and $\theta_{d\over2}$ is the $d\over2$-dimensional null column-vector. Since $\left(\sum_{i=1}^nM^K_n\dots M^K_{i+1}M^{K,K^c}_iM^{K^c}_{i-1}\dots M^{K^c}_1\right)W=\left(\sum_{i=1}^nd_{i-1}\dots d_1\right)W$ diverges, $M_n\dots M_1$ also do although $M^{K^c}_n\dots M^{K^c}_1$ converges to the null matrix.

\end{rem}

\begin{rem}The two conditions of (ii) are satisfied -- assuming that the ones of (i) are -- if all the submatrices $M_n^{K^c}$ have spectral radius less than $1$, or if their eigenvalues greater or equal to $1$ disappear in the product $M^{K^c}_{n+h}\dots M^{K^c}_n$ for some fixed $h$ and for any $n$.

\end{rem}

\begin{rem}Let us compare now the problem of the convergence of $Q_n=M_n\dots M_1$ to the one of
$$
R_n:={M_n\dots M_1\over\Vert M_n\dots M_1\Vert}.
$$
Let $(M_n)$ be a sequence of complex-valued matrices such that $R_n$ converges. Since the limit matrix $R$ has norm $1$, it has some nonnull columns; let us prove that they are right-eigenvectors of each matrix $M$ that occurs infinitely many times in the sequence $(M_n)$. The nonnegative real $\displaystyle\lambda_n:=\Vert MR_{n-1}\Vert$ is bounded by $\Vert M\Vert$ and satisfy $MR_{n-1}=\lambda_nR_n$ for any $n$ such that $M_n=M$, so it has at least one limit point $\lambda$ that satisfy $MR=\lambda R$, and the columns of $R$ are right-eigenvectors of $M$ for the eigenvalue $\lambda$.

Suppose that $\lambda\ne0$ -- perhaps it is not possible that $\lambda=0$. The convergence of $R_n$ can hold without the convergence of $\displaystyle\left({1\over\lambda}M_n\right)\dots\left({1\over\lambda}M_1\right)$, see for instance the case where all the $M_n$ are equal to $\left(\begin{array}{cc}1&1\\0&1\end{array}\right)$. But conversely the convergence of this last product to a nonnull matrix, for some $\lambda\ne0$, implies obviously the one of $R_n$.

\end{rem}




\begin{thebibliography}{99}

\bibitem{DL} \textsc{I. Daubechies \& J. C. Lagarias},
Sets of matrices all infinite products of which converge,
{\sl Linear Algebra and its Applications} {\bf 161} (1992), 227-263.
\bibitem{DL'} \textsc{I. Daubechies \& J. C. Lagarias},
Corrigendum/addendum to: Sets of matrices all infinite products of which
converge,
{\sl Linear Algebra and its Applications} {\bf 327} (2001), 69-83.
\bibitem{EF} \textsc{L. Elsner \& S. Friedland},
Norm conditions for convergence of infinite products,
{\sl Linear Algebra and its Applications} {\bf 250} (1997), 133-142, {\href{http://www.math.uni-bielefeld.de/~elsner/elsner-publications.html}{\sl download the list of publications}}
\bibitem{Erd} \textsc{P. Erd\H os},
On a family of symmetric Bernoulli convolutions.
{\sl Amer. J. Math.} {\bf61} (1939), 974-976, \href{http://www.renyi.hu/~p_erdos/1939-05.pdf}{\sl download}.
\bibitem{L} \textsc{J. Lorenz},
Convergence of products of stochastic matrices with positive diagonal and the opinion dynamics background (2007), {\sl\href{http://arxiv.org/pdf/0708.3177v1}{download}}.
\bibitem{MNR} \textsc{A. Mukherjea, A. Nakassis \& J.S. Ratti},
On the distribution of the limit of products of i. i. d. $2\times2$ random stochastic matrices, 
{\sl J. Theor. Probab.} {\bf 12} (1999), 571-583.
\bibitem{MN} \textsc{A. Mukherjea \& A. Nakassis},
On the continuous singularity of the limit distribution of products of i. i. d. $d\times d$ 
stochastic matrices,
{\sl J. Theor. Probab.} {\bf 15} (2002), 903-918.
\bibitem{OST} \textsc{E. Olivier, N. Sidorov, \& A. Thomas},
On the Gibbs properties of Bernoulli convolutions
related to $\beta$-numeration in multinacci bases,
{\sl Monatsh. Math.} {\bf 145} (2005), 145--174, \href{http://www.cmi.univ-mrs.fr/~thomas/publications/OST_Monat-F-Math.pdf}{\sl download}.
\bibitem{S} \textsc{E. Seneta},
Non-negative matrices and Markov chains,
{\sl Springer Series in Statistics. New York - Heidelberg - Berlin: Springer-
Verlag} {\bf XV} (1981), \href{http://books.google.fr/books?id=hsE0__8frPoC&printsec=frontcover&dq=seneta+nonnegative&cd=1#v=onepage&q=seneta%20nonnegative&f=false}{\sl partial download}.

\end{thebibliography}
\end{document}